\documentclass[a4paper,notitlepage, 12pt,reqno]{article}

  \usepackage{mathtext}
 \usepackage[cp1251]{inputenc}
  \usepackage[T2A]{fontenc}
 \usepackage[russian,english]{babel}
 \usepackage[all,arc,poly,2cell,curve,arrow,tips]{xy}

  \usepackage{enumerate, eucal, amsthm,amsmath, amssymb}

  \allowdisplaybreaks[4]

 \theoremstyle{definition}

 \theoremstyle{plain}
 \newtheorem{thm}{Theorem}
 \newtheorem*{thm*}{Theorem}
 \newtheorem{prop}{Proposition}
  \newtheorem*{prop*}{Предложение}
 \newtheorem{cor}{Corrolary}
  \newtheorem*{cor*}{Следствие}
 
  \newtheorem*{lem*}{Лемма}

 \theoremstyle{remark}

 \newtheorem*{remark*}{Замечание}

 \renewcommand{\abstractname}{}

  \newcounter{ab}
\title{ The Gelfand-Tsetlin-Zhelobenko  base vectors for the series $B$}

 \author{D.V. Artamonov\footnote{Lomonosov Moscow State University, artamonov.dmitri@gmail.com}}
  \date{}

\begin{document}
 \maketitle

\renewcommand{\abstractname}{}

\begin{abstract}

 Using the method of $Z$-invariants of  Zhelobenko we construct  base vectors  of Gelfand-Tsetlin type in  the space of $\mathfrak{o}_{2n-1}$-highest vectors in a representation of $\mathfrak{o}_{2n+1}$.
The construction is based on a  relation between restriction problems  $\mathfrak{o}_{2n+1}\downarrow\mathfrak{o}_{2n-1}$ and  $\mathfrak{gl}_{n+1}\downarrow\mathfrak{gl}_{n-1}$.

\end{abstract}

\section{Introduction}
In the paper \cite{GC2} Gelfand and Tsetlin constructed a base in a representation of the Lie algebra  $\mathfrak{o}_N$.  The construction is based on an investigation of a branching of an irrep of  $\mathfrak{o}_N$ under the restriction of Lie algebras $\mathfrak{o}_N\downarrow\mathfrak{o}_{N-1}$. The restriction problem $\mathfrak{g}\downarrow\mathfrak{k}$,  where $\mathfrak{k}$ is a subalgebra in a Lie algebra $\mathfrak{g}$, is a problem of an explicit description of   $\mathfrak{k}$-highest vectors in a representation of $\mathfrak{g}$.

 Later it turned out that it is natural to have a construction of a base of Gelfand-Tselin type for a representation of $\mathfrak{o}_{2n+1}$,  based on restrictions $\mathfrak{o}_{2n+1}\downarrow\mathfrak{o}_{2n-1}$ inside the series $B$.
Thus in  physical literature there were attempts to obtain a such a base for
$\mathfrak{o}_5$. Such a base is needed in the problem of classification of states of a five-dimensional quasi-spin in a  shell models of nuclear kernels  \cite{h}, \cite{hl}, \cite{ag}.

 Zhelobenko in \cite{zh2}  constructed base vectors \footnote{ But he did not manage to obtain  formulas for the action of generators of the algebra in this base} of Gelfand-Tsetlin type for $\mathfrak{sp}_{2n}$. He used a simpler technique of  $Z$-invariants.  This technique allowed him  to find a relation between the restriction problems $\mathfrak{sp}_{2n}\downarrow\mathfrak{sp}_{2n-2}$ and $\mathfrak{gl}_{n+1}\downarrow\mathfrak{gl}_{n-1}$.


Later V.V. Shtepin investigated the problem of restriction  $\mathfrak{o}_{2n+1}\downarrow\mathfrak{o}_{2n-1}$ in  \cite{sh2} using the technique of  $Z$-invariants of Zhelobenko, but he did not find a relation with the problem of restriction $\mathfrak{gl}_{n+1}\downarrow\mathfrak{gl}_{n-1}$.

Finally the problem of construction of a Gelfand-Tsetlin type base for the series $B$,  $D$, and $C$  was solved completely by Molev (see \cite{M}). To obtain such a construction a solution of restriction problems were obtained. Molev constructed base vectors and obtained formulas for the action of generators of the algebras in this base.  But he used a much more difficult technique. The key step  in the Molev's approach is a construction of an action of a Yangian on the space of  $\mathfrak{o}_{2n-1}$-highest vectors with a fixed highest weight.  Also he did not point out a relation of restriction problems $\mathfrak{o}_{2n+1}\downarrow\mathfrak{o}_{2n-1}$  and  $\mathfrak{gl}_{n+1}\downarrow\mathfrak{gl}_{n-1}$.

In the present paper in Section \ref{gczs0}  using the technique of $Z$-invariants of  Zhelobenko we find  a relation between the problems of restriction $\mathfrak{o}_{2n+1}\downarrow \mathfrak{o}_{2n-1}$ and $\mathfrak{gl}_{n+1}\downarrow\mathfrak{gl}_{n-1}$. Using it in Section \ref{gczs1} we construct in the same manner as in the case $\mathfrak{sp}_{2n}$ in \cite{zh2} base vectors of Gelfand-Tsetlin type for the algebra  $\mathfrak{o}_{2n+1}$ (Theorem \ref{osnt}).
   Unfortunately we have not managed to find the formulas for the action of generators of the algebra.
    The structure of our  Gelfand-Tsetlin tableaux for $\mathfrak{o}_{2n+1}$ is the same as the structure of  Gelfand-Tsetlin tableaux constructed by Molev (see \cite{M}).

\section{
The algebra  $\mathfrak{o}_{2n+1}$, the method of $Z$-invariants}


The algebra $\mathfrak{o}_{2n+1}$ is generated by  $(2n+1)\times (2n+1)$ matrices, whose rows and columns are indexed by  $-n,...,-1,0,1,...,n$, of type
\begin{align*}&F_{i,j}=E_{i,j}-E_{-j,-i},\,\,\,\,\,\, i,j=-n,...,-1,0,1,...,n,\end{align*}

where $E_{i,j}$ is a matrix unit.

The subalgebra
$\mathfrak{o}_{2n-1}$ is generated by $F_{i,j}$
for $i,j\in -n,...,-2,0,2,...,n$.

\section{Zhelobenko's realization }


On the open dense subset $O_{2n+1}^0\subset O_{2n+1}$ the Gauss
decomposition takes place

\begin{align*}&O_{2n+1}^0=Z^{-}DZ,\,\,\, X=\zeta \delta z,\,\,\,\\& X\in O_{2n+1},\,\,\, \zeta\in Z^{-},\,\,\, \delta\in D,\,\,\, z\in
Z,\end{align*}

where $Z^{-}$ is a subgroup of upper-triangular unipotent matrices
form $O_{2n+1}$, $D$ is a subgroup of diagonal matrices in
$O_{2n+1}$ , and $Z$ is a subgroup of lower-triangular unipotent
matrices. On the space of polynomial functions on
   $Z$ there exists an action of  $O_{2n+1}$  by the following ruler.
   Let us be given a function on  $Z$ of type $f(z)=f(z_{i,j})$, $i<j$. For  $X\in   O_{2n+1}$
   put

\begin{equation}\label{xfz}(Xf)(z)=\alpha(\widetilde{\delta})f(\widetilde{z}),\,\,\,
zX=\widetilde{\zeta}\widetilde{\delta}\widetilde{z},\,\,\,\,\,\alpha(\delta)=\delta_{-n}^{r_{-n}}...\delta_{-1}^{r_{-1}},\end{equation}

where $\delta=diag(\delta_{-n},\delta_{-n+1},...,\delta_{n})$

   Thus the space of all such  functions form a representation of $O_{2n+1}$.  A finite-dimensional representation with the highest weight
  $[m_{-n},...,m_{-1}]$, where numbers  $m_{-i}$ are simultaneously integers or half-integers,  is formed by functions that satisfy a system of PDE called   {\it  the indicator system}:
$$L_{-n,-n+1}^{r_{-n}+1}f=0,...,L_{-1,0}^{r_{-1}+1}f=0,$$

where $r_i$ are defined as follows

\begin{equation}\label{rb}r_{-n}=m_{-n}-m_{-n+1},...,r_{-2}=m_{-2}-m_{-1},\,\, r_{-1}=2m_{-1}.\end{equation}

 Here  $L_{i,j}$ are operator that do
left infinitesimal shifts of a function $f(z)$ by $F_{i,j}$.

The procedure of a construction of the Gelfand-Tsetlin type base is
based  on an investigation  of a branching of an irrep under the
restriction of the algebra. The method of $Z$-invariants gives us a description of functions that are $\mathfrak{o}_{2n-1}$-highest vectors.
As in   the case  $\mathfrak{sp}_{2n}$ (see  \cite{zh2})  one can easily show that   a function  $f$ is a $\mathfrak{o}_{2n+1}$-highest vector $f$  if and only if the following conditions hold.
\begin{enumerate}
\item The function $f$ depends on the following variables \begin{equation}\label{zb}f=f(z_{-n,-1},...,z_{-2,-1},z_{-n,1},...,z_{-2,1},z_{0,1}).\end{equation} We used the relation $z_{-1,1}=-\frac{z_{0,1}^2}{2},$ that holds for
the matrix elements of the group $Z$.
\item The function $f$ satisfies the indicator system.
 \end{enumerate}


\section{An explicit form of the indicator system and it's solutions}

Let us write the explicit form of the indicator system being restricted to the functions of type  \eqref{zb}.
The indicator system looks as follows

\begin{align}\begin{split}\label{indicb}
&L_{-n,-n+1}^{r_{-n}+1}f=(z_{-n+1,-1}\frac{\partial}{\partial z_{-n,-1
}}+z_{-n+1,1}\frac{\partial}{\partial
z_{-n,1}})^{r_{-n}+1}f=0,\\&...\\
&L_{-3,-2}^{r_{-3}+1}f=(z_{-2,-1}\frac{\partial}{\partial  z_{-3,-1}}+z_{-2,1}\frac{\partial}{\partial z_{-3,1}})^{r_{-3}+1}f=0\\
&L_{-2,-1}^{r_{-2}+1}f=(\frac{\partial}{\partial  z_{-2,-1}}+\frac{z_{0,1}^2}{2}\frac{\partial}{\partial z_{-2,1}})^{r_{-2}+1}f=0\\
&L_{-1,0}^{r_{-1}+1}f=(\frac{\partial}{\partial z_{0,1}})^{r_{-1}+1}f=0.\end{split}
\end{align}

To solve it let us introduce new variables \begin{align*}&
u_{-k}=z_{-k,1}+\frac{z_{0,1}^2}{2}z_{-k,-1},\,\,\,v_k=z_{-k,1}-\frac{z_{0,1}^2}{2}z_{-k,-1},\,\,k=2,...,n.\\&
u_{-1}=z_{0,1},\,\,\, v_{-1}=0.
\end{align*}

The variables
$z_{-k,-1},z_{-k,1}$ and be reconstructed as follows $u_{-k},v_{-k}$:

\begin{equation*}
z_{-k,1}=\frac{u_{-k}+v_{-k}}{2},\,\,\,z_{-k,-1}=\frac{u_{-k}-v_{-k}}{\frac{z_{0,1}^2}{2}}.
\end{equation*}

In the space of polynomials in variables  $z_{-k,-1},z_{-k,1},...,z_{-2,-1},z_{-2,1},z_{0,1}$ there exists a base

\begin{align}
\begin{split}
\label{polyn}
&u_{-1}^{p_{-1}}\prod_{k=2}^{n}(u_{-k}+v_{-k})^{p_{-k}}(u_{-k}-v_{-k})^{q_{-k}},\\
&p_{-k},q_{-k}\geq 0,\,\,\,\,\,k=2,...,n,    \,\,\,\,\,p_{-1}+2\sum_{k=2}^np_{-k}\geq 0.
\end{split}
\end{align}

Let us find a condition under which   \eqref{polyn} is a solution.

Consider first the equations $L_{-k,-k+1}^{r_{-k}+1}f=0$ for  $k=n,...,2$.
One has
\begin{align*}
&L_{-k,-k+1}u_{-k}=u_{-k+1},\,\,\, L_{-k,-k+1}v_{-k}=v_{-k+1},\\& L_{-k,-k+1}u_{-l}=L_{-k,-k+1}v_{-l}=0,\,\,\,
k\neq l,\,\, k=2,..,n,\end{align*}

Thus the operator  $L_{-k,-k+1}^{r_{-k}+1}$  maps a polynomial in  variables $u_{-i},v_{-i}$ into zero if and only if in each monomial the sum of degrees of $u_{-k} $  and  $v_{-k}$ is not greater than  $r_{-k}$. That is if

 $$p_{-k}+q_{-k}\leq r_{-k},\,\,\,\,\,k=2,...,n.  $$

Consider the equation $L_{-1,0}^{r_{-1}+1}f=0$. One has
\begin{align*}
 &L_{-1,0}=\frac{\partial}{\partial z_{0,1}},
\end{align*}

the operator  $L_{-1,0}^{r_{-1}+1}$ maps a polynomial in variables $u_{-i},v_{-i}$  into zero if and only if  in variables $z$ we have a polynomial in variable $z_{0,1}$ of degree not greater than $r_{-1}$.  The polynomial  \eqref{polyn} being rewritten  in variables  $z$ has a degree in variable $z_{0,1}$ equal to  $p_{-1}+2\sum_{k=2}^{n}(p_{-k}+q_{-k})$. That is the following condition must hold

$$
p_{-1}+2\sum_{k=2}^{n}(p_{-k}+q_{-k})\leq r_{-1}
$$

Thus we obtain

\begin{prop}
In the space of solutions of \eqref{indicb} there exists a base

\begin{equation}\label{solindib0}
f=u_{-1}^{p_{-1}}\prod_{k=2}^{n}(u_{-k}+v_{-k})^{p_{-k}}(u_{-k}-v_{-k})^{q_{-k}},
\end{equation}

where \begin{align}
\begin{split}
\label{solindib1}
&p_{-k},q_{-k}\geq 0,\,\,\,\,\,k=2,...,n,    \,\,\,\,\,p_{-1}+2\sum_{k=2}^np_{-k}\geq 0,\\
& p_{-k}+q_{-k}\leq
r_{-k},\,\,k=2,...,n,\,\,\,
p_{-1}+2\sum_{k=2}^{n}(p_{-k}+q_{-k})\leq r_{-1}.
\end{split}
\end{align}

\end{prop}

Let us find the action of $F_{-i,-i}$ onto these functions. The matrix  $\widetilde{z}$ from \eqref{xfz} for  $X=e^{tF_{-i,-i}}$ can be obtained from  $z$ by multiplication of the row  $-i$  onto  $e^{-t}$, of the row  $i$ onto  $e^{t}$,  of the column  $-i$ onto  $e^{t}$, of the column  $i$ onto $e^{-t}$. The matrix $\delta$ equals $e^{tF_{-i,-i}}$. Thus for the infinitesimal action one has

\begin{align}
\begin{split}
\label{fii}
&F_{-i,-i}f=-z_{-i,-1}\frac{\partial}{\partial z_{-i,-1}}f-z_{-i,1}\frac{\partial}{\partial z_{-i,1}}f+m_{-i}f,\,\,\,\, i=2,...,n,\\
&F_{-1,-1}f=\sum_{i=2}^n(z_{-i,-1}\frac{\partial}{\partial z_{-i,-1}}-z_{-i,1}\frac{\partial}{\partial z_{-i,1}})f-z_{0,1}\frac{\partial}{\partial z_{0,1}}f+m_{-1}f
\end{split}
\end{align}

\subsection{The restriction problem $\mathfrak{gl}_{n+1}\downarrow\mathfrak{gl}_{n-1}$}

\subsubsection{The indicator system and it's solutions}

Consider the algebra of all matrices  $\mathfrak{gl}_{n+1}$ acting in the space with coordinates indexed by  $-n,...,-1,1$. Representations of this Lie algebra can be realized in the space of functions on upper-triangular unipotent matrices. An irreducible representation is selected by an indicator system (see \cite{zh2}).

Consider the subalgebra  $\mathfrak{gl}_{n-1}$ generated by $E_{i,j}$,  $i,j\in\{-n,...,-2\}$. Let us b given an irreducible representation of   $\mathfrak{gl}_{n+1}$ with the highest weight $[m_{-n},...,m_{-1},m_{1}=0]$. Consider the problem of restriction  $\mathfrak{gl}_{n+1}\downarrow\mathfrak{gl}_{n-1}$.

It turns out that (see \cite{zh2}) $\mathfrak{gl}_{n-1}$-highest vectors are functions of type $$f=f(z_{-n,-1},...,z_{-2,-1},z_{-n,1},...,z_{-1,1}),$$ satisfying the indicator system.

Being restricted to these function the indicator system  takes the following explicit form

\begin{align}\begin{split}\label{indicg}
&L_{-n,-n+1}^{r_{-n}+1}f=(z_{-n+1,-1}\frac{\partial}{\partial z_{-n,-1
}}+z_{-n+1,1}\frac{\partial}{\partial
z_{-n,1}})^{r_{-n}+1}f=0,\\&...\\
&L_{-3,-2}^{r_{-3}+1}f=(z_{-2,-1}\frac{\partial}{\partial  z_{-3,-1}}+z_{-2,1}\frac{\partial}{\partial z_{-3,1}})^{r_{-3}+1}f=0\\
&L_{-2,-1}^{r_{-2}+1}f=(\frac{\partial}{\partial  z_{-2,-1}}+z_{-1,1}\frac{\partial}{\partial z_{-2,1}})^{r_{-2}+1}f=0\\
&L_{-1,1}^{r_{-1}+1}f=(\frac{\partial}{\partial z_{-1,1}})^{r_{-1}+1}f=0,\end{split}
\end{align}

where
\begin{equation}\label{ra}r_{-n}=m_{-n}-m_{-n+1},...,r_{-2}=m_{-2}-m_{-1},\,\, r_{-1}=m_{-1}.\end{equation}

To solve this system let us introduce new variables

 \begin{align*}
&x_{-k}=z_{-k,1}+z_{-1,1}z_{-k,-1},\,\,\,y_{-k}=z_{-k,1}-z_{-1,1}z_{-k,-1},\,\,k=2,...,n.\\&
x_{-1}=z_{-1,1},\,\,\, y_{-1}=0.
\end{align*}

The variables
$z_{-k,-1},z_{-k,1}$ can be reconstructed from  $x_{-k},
y_{-k}$:

\begin{equation}
z_{-k,1}=\frac{x_{-k}+y_{-k}}{2},\,\,\,z_{-k,-1}=\frac{x_{-k}-y_{-k}}{2z_{-1,1}}.
\end{equation}

Thus the space of polynomials   $z_{-k,-1},z_{-k,1},...,z_{-2,-1},z_{-2,1},z_{-1,1}$ there exists a base

\begin{align}
\begin{split}
\label{polyn}
&x_{-1}^{p_{-1}}\prod_{k=2}^{n}(x_{-k}+y_{-k})^{p_{-k}}(x_{-k}-y_{-k})^{q_{-k}},\\
&p_{-k},q_{-k}\geq 0,\,\,\,\,\,k=2,...,n,    \,\,\,\,\,p_{-1}+\sum_{k=2}^np_{-k}\geq 0.
\end{split}
\end{align}

Let us write conditions under which this polynomial is a solution.
One has

\begin{align*}
&L_{-k,-k+1}x_{-k}=x_{-k+1},\,\,\,\ L_{-k}y_{-k}=y_{-k+1},\,\,\,  \\&L_{-k,-k+1}x_{-l}=L_{-k,-k+1}y_{-l}=0,\,\,\,
k\neq l,\,\, k=2,..,n.
\end{align*}

As  in the previous Section we obtain the following statement.

\begin{prop} In the space of polynomial solutions of the system \eqref{indicg}
there exists a base of type

\begin{equation}\label{solg0}
f=x_{-1}^{p_{-1}}\prod_{k=2}^{n}(x_{-k}+y_{-k})^{p_{-k}}(x_{-k}-y_{-k})^{q_{-k}},
\end{equation}

where \begin{align}
\begin{split}
\label{solg1}
&p_{-k},q_{-k}\geq 0,\,\,\,\,\,k=2,...,n,    \,\,\,\,\,p_{-1}+\sum_{k=2}^np_{-k}\geq 0,\\
& p_{-k}+q_{-k}\leq r_{-k},\,\,k=2,...,n,\,\,\, p_{-1}+\sum_{k=2}^{n}(p_{-k}+q_{-k})\leq r_{-1}.
 \end{split}
\end{align}

\end{prop}

Let us find the action of $E_{-i,-i}$ onto these functions. The matrix $\widetilde{z}$ from \eqref{xfz} for  $X=e^{tE_{-i,-i}}$ can be obtained form   $z$ by multiplication of the row  $-i$  onto  $e^{-t}$ and by multiplication of the column  $-i$ onto  $e^{t}$. The matrix $\delta$ equals  $e^{tE_{-i,-i}}$. Thus for an infinitesimal action one has

\begin{align}
\begin{split}
\label{eii}
&E_{-i,-i}f=-z_{-i,-1}\frac{\partial}{\partial z_{-i,-1}}f-z_{-i,1}\frac{\partial}{\partial z_{-i,1}}f+m_{-i}f,\,\,\,\, i=2,...,n,\\
&E_{-1,-1}f=\sum_{i=2}^nz_{-i,-1}\frac{\partial}{\partial z_{-i,-1}}
f-z_{-1,1}\frac{\partial}{\partial z_{-1,1}}f+m_{-1}f,\\
&E_{1,1}f=\sum_{i=1}^nz_{-i,1}\frac{\partial}{\partial z_{-i,1}}
f.
\end{split}
\end{align}

\subsubsection{The Gelfand-Tsetlin base.}

In the space of $\mathfrak{gl}_{n-1}$-highest vectors there exists
the Gelfand-Tsetlin base  encoded by  tableaux in which the
betweeness conditions hold.

\begin{displaymath}
 \xymatrix{ m_{-n}   && m_{-n+1} &  ... & m_{-1} && 0\\
 &m'_{-n,n}   && ...   && m'_{-n,-1}\\
&& m_{-n,n-1} &...& m_{-2,n-1}\\
}
\end{displaymath}

To prove the main statement below we need a realization of a
representation in the functions on the whole group. Onto a function
$f(g)$ an element $X\in GL_{n}$  acts by the ruler

$$(Xf)(g)=f(gX).$$

Let  $a^{j}_{i}$ be a function of a matrix element, where $j$ is a row index and $i$ is a column index. Put
$$a_{i_1,...,i_k}:=det(a_i^j)_{i=i_1,...,i_k}^{j=-n,...,-n+k-1}$$
One can easily check that the function
\begin{equation}
v_0=\label{stv} \prod_{k=-n}^{-1} (a_{-n,...,-k})^{r_{-k}},
\end{equation}
 is a highest vector for $\mathfrak{gl}_{n+1}$ with the weight  $[m_{-n},...,m_{-1},0]$.
  Indeed the operator  $E_{i,j}$acts onto  $a_{i_1,...,i_k}$ by the ruler

$$
a_{i_1,...,i_k}\mapsto a_{\{i_1,...,i_k\}\mid_{j\mapsto i}},
$$
where $.\mid_{j\mapsto i}$ is an operation of substitution of the index $i$ instead of  $j$,  if $j\notin\{i_1,...,i_k\}$ then the determinant is mapped to zero. Onto a product of determinant the operator $E_{i,j}$ act by the Leibnitz ruler.

To write the formulas for a vector corresponding to a tableau let us
introduce operators $e_{1,-i}$, $i=n,...,1$ acting onto determinants
by the ruler

\begin{align*}
&a_{-n,...,-i-1,-i}\mapsto a_{-n,...,-i-1,1},\\
&\text{other determinants}\mapsto 0,
\end{align*}
and acting onto a product of determinants by the Leibnitz ruler.
Also let us introduce operators $e_{-1,-i}$, $i=n,...,2$ acting onto
determinant by the ruler

\begin{align*}
&a_{-n,...,-i-1,-i}\mapsto a_{-n,...,-i-1,-1},\\
&a_{-n,...,-i-1,-i,1}\mapsto a_{-n,...,-i-1,-1,1},\\
&\text{other determinants}\mapsto 0,
\end{align*}
and acting onto a product of determinants by the Leibnitz ruler.

Then for a  $\mathfrak{gl}_{n-1}$-highest vector $v$, defined by a Gelfand-Tsetlin
tableau  one has a formula

\begin{equation}
\label{fsv}
v=const\cdot\prod_{i=-n}^{-2}e_{-1,-i}^{m'_{-i,n}-m_{-i,n-1}}\prod_{i=-n}^{-1}e_{1,-i}^{m_{-i,n}-m'_{-i,n}}v_0.
\end{equation}

See for example \cite{m2}. Our operators  $e_{\pm 1,-i}$ correspond
to operators $pE_{\pm 1,-i}$\footnote{Note that in \cite{m2} and in
the present paper the indexation of coordinates is different}.

Indeed  the extremal projector  $p$ maps a vector  $v$ to zero in the case  $v=E_{-}w$, where $E_{-}$ is an element of  $\mathfrak{gl}_{n-1}$, corresponding to a negative root.
The coincidence of actions of  $pE_{\pm 1,-i}$ and  $e_{\pm 1,-i}$ onto determinant can be easily checked. One has to prove that their action onto products of determinants coincide. Let us write it as follows
\begin{equation}
\label{proi}
\prod a_{-n,...,-k}^{\alpha_k}a_{-n,...,-k,-1}^{\beta_k}a_{-n,...,-k,1}^{\gamma_k}a_{-n,...,-k,-1,1}^{\delta_k}
\end{equation}

The operator  $E_{ -1,-2}$  act only on factors with  $k=2$, the application of   $p$ changes nothing. The resulting action of  $pE_{ -1,-2}$  coincides with $e_{-1,-2}$. Now consider
$E_{ -1,-3}$, this operator can be represented as  $[E_{-1,-2},E_{-2,-3}]=E_{-1,-2}E_{-2,-3}-E_{-2,-3}E_{-1,-2}$.  After application of $p$ we obtain  $pE_{-1,-2}E_{-2,-3}$. The operator  $E_{-2,-3}$ act onto determinants in \eqref{proi} with  $k=3$. Under the action of  $E_{ -1,-2}E_{-2,-3}$ one obtains a product of determinant of type \eqref{proi}, which consists of determinants that are highest with respect to  $\mathfrak{gl}_{n-1}$. Thus   $pE_{ -1,-3}=E_{ -1,-2}E_{-2,-3}$. One can easily prove that  $E_{ -1,-2}E_{-2,-3}$ equals to $e_{-1,-3}$. Thus finally $pE_{-1,-3}=e_{-1,-3}$.
For the rest operators $e_{\pm 1,-i}$ the proof is the same.

 In
the formula  analogous to \eqref{fsv} not $pE_{\pm 1,-i}$ but the operators
denoted in \cite{m2} as $z_{\pm 1,-i}$ occur. However under the
action onto weight vectors the operators $pE_{\pm 1,-i}$ and $z_{\pm
1,-i}$  are proportional. Thus the formula \eqref{fsv} follows from
the results of \cite{m2}.

We need the following statement.

\begin{prop}
\label{p1}
If one decomposes a vector corresponding to  a Gelfand-Tsetlin tableau
 by the base \eqref{solg0}, then for the summands one
has the equality $p_{-1}+\sum_{k=-2}^n(p_{-k}+q_{-k})\leq
m_{-1,n}-m'_{-1,n}$. For at least one of the summand the equality takes place.
\end{prop}
\proof

Take a realization of a representation on the functions on the whole group.
From explicit formulas for the action of $e_{\pm 1,-i}$ one obtains that
$v$ defined by \eqref{fsv} is a linear combination of products of determinants of type
$a_{-n,...,-i}$, $a_{-n,...,-i,1}$, $a_{-n,...,-i,-1}$,
$a_{-n,...,-i,-1,1}$. In which such a product the sum of degrees of
$a_{-n,...,-2,1}$ and $a_{-n,...,-3,-1,1}$ equals to
$m_{-1,n}-m'_{-1,n}$.

Indeed these determinants appear as a result of action of  $e_{1,-1}$ and  $e_{-1,-2}$. The sum of degrees of theese operatoes in \eqref{fsv} equals  $m_{-1,n}-m'_{-1,n}$.
Thus the sum of degrees of $a_{-n,...,-2,1}$  and  $a_{-n,...,-3,-1,1}$ in  \eqref{fsv} equals  $m_{-1,n}-m'_{-1,n}$.

Consider the value of these determinants on the subgroup $Z$. This is
a polynomial in variables  $z_{-n,1},z_{-n,-1}$, and $z_{-1,1}$ can
appear only from the determinants $a_{-n,...,-2,1}$ and
$a_{-n,...,-3,-1,1}$. More presice $a_{-n,...,-2,1}\mid_Z=z_{-1,1}$,  $a_{-n,...,-3,-1,1}\mid_{Z}=z_{-2,-1}z_{-1,1}-z_{-2,1}$.  Thus it's degree in the variable  $z_{-1,1}$
equals to $m_{-1,n}-m'_{-1,n}$.

Now let us pass to the variables  $x_{-i}$, $y_{-i}$. Let us be
given a base polynomial \eqref{solg0} in these variables, rewrite it in variables
$z_{-n,1},z_{-n,-1}$, then it's degree in $z_{-1,1}$ equals to
$p_{-1}+\sum_{k=-2}^n(p_{-k}+q_{-k})$.

 Thus when we pass from   $z_{-n,1},z_{-n,-1}$ to  $x_{-i}$, $y_{-i}$
 only monomials with   $p_{-1}+\sum_{k=-2}^n(p_{-k}+q_{-k})\leq m_{-1,n}-m'_{-1,n}$ appear and at least for one of them the equality takes
 place.

\endproof

\begin{cor}
The span of vectors \eqref{solg0} with  $p_{-1}+\sum_{k=-2}^n(p_{-k}+q_{-k})\leq m_{-1,n}-1$ contains a subspace the span of vectors corresponding to tableaux with  $m'_{-1,n}>0$.
\end{cor}

\section{A relation between restriction problems}
\label{gczs0}

Let us establish a relation between restriction problems   $\mathfrak{o}_{2n+1}\downarrow
\mathfrak{o}_{2n-1}$ with exponents
\eqref{rb} and $\mathfrak{gl}_{n+1}\downarrow\mathfrak{gl}_{n-1}$ with exponents
\eqref{ra}.




Consider the cases when $p_{-1}$ is odd and even separately.

\subsection{The case of even $p_{-1}$} Let us write $p_{-1}=2p'_{-1}$.
Take the solution \eqref{solindib0} with exponents
$r_{-n},...,r_{-2},r_{-1}$ and relate to it a solution \eqref{solg0}
with exponents $r_{-n},...,r_{-2},[\frac{r_{-1}}{2}]$, where $[.]$ is an integer part, according to the
ruler

\begin{equation}
\label{sop1} (p_{-1},
p_{k},q_{k}),\,\,\, k=-2,...,-n \mapsto  p'_{-1}=\frac{p_{-1}}{2},
p_{k},q_{k},\,\,\, k=-2,...,-n
\end{equation}

The inequalities \eqref{solindib1} for $p_{-1}$, $p_{k},q_{k}$, and
$r_{-n},...,r_{-2},r_{-1}$ are equivalent to inequalities
\eqref{solg1} with $p'_{-1}$, $p_{k},q_{k}$ and
$r_{-n},...,r_{-2},[\frac{r_{-1}}{2}]$. Thus we obtain the following
statement.

\begin{prop}
The correspondence  \eqref{sop1}  is a bijection between the
solution space \eqref{indicb} with even $p_{-1}$ and the space of
all solutions \eqref{indicg}.
\end{prop}


\subsection{The case of odd $p_{-1}$} Let us write $p_{-1}=2p'_{-1}+1$.
Take a solution \eqref{solindib0} of the system with exponents
$r_{-n},...,r_{-2},r_{-1}$ and let us relate to it a solution
\eqref{solg0} with exponents $r_{-n},...,r_{-2},[\frac{r_{-1}}{2}]$
according to the ruler

\begin{equation}
\label{sop2} (p_{-1},
p_{k},q_{k}),\,\,\, k=-2,...,-n \mapsto  (p'_{-1}=\frac{p_{-1}-1}{2},
p_{k},q_{k}),\,\,\, k=-2,...,-n
\end{equation}
Let us first prove that this correspondence is well defined. It is necessary to check that for the image of  \eqref{sop2} the inequalities  \eqref{solg1} hold.
Inequalities \eqref{solindib1} for $p_{-1}$, $p_{k},q_{k}$ and
$r_{-n},...,r_{-2},r_{-1}$ give us inequalities \eqref{solg1} for
$p'_{-1}$, $p_{k},q_{k}$ and $r_{-n},...,r_{-2},[\frac{r_{-1}}{2}]$.
That is \eqref{sop2} is a well-defined embedding of the space of
solutions with even   $p_{-1}$ of the  problem of restriction
$\mathfrak{o}_{2n+1}\downarrow \mathfrak{o}_{2n-1}$ into the space
of solutions of the problem of restriction
$\mathfrak{gl}_{n+1}\downarrow \mathfrak{gl}_{n-1}$.

 Let us
describe the image of this embedding.

\begin{prop}
If the highest weight of the representation of $\mathfrak{o}_{2n+1}$
is half-integer then \eqref{sop2} is an isomorphism.

If the highest weight of the representation of $\mathfrak{o}_{2n+1}$
is integer then the image of \eqref{sop2} is a span of tableaux with
$m'_{-1,n}>0$.
\end{prop}

\proof

Take the inequality  $p_{-1}+2\sum_{k=-2}^{-n}(p_k+q_k)\leq r_{-1},$
and divide it by two, one obtains
\begin{equation}\label{kl}0\leq p'_{-1}+\frac{1}{2}+\sum_{k=-2}^{-n}(p_k+q_k)\leq \frac{r_{-1}}{2}.\end{equation}

Suggest that the highest weight is half-integer, that is  $r_{-1}$ is odd.
Because of the fact that $p'_{-1}$, $p_{-i}$, $q_{-i}$ are integer,
the inequality \eqref{kl} is equivalent to the following one:
$p'_{-1}+\sum_{k=-2}^{-n}(p_k+q_k)\leq \frac{r_{-1}-1}{2}=[\frac{r_{-1}}{2}].$  Thus from
the equality \eqref{solg1} the equality
 \eqref{solindib1} follows. That is \eqref{sop2} is an isomorphism.

Suggest that the highest weight is integer, that is
$r_{-1}$ is even.

Let us first prove that the image of the correspondence contains the linear span of tableaux with $m'_{-1,n}>0$.
We obtain the equality
\begin{equation}\label{ner}p'_{-1}+\frac{1}{2}+\sum_{k=-2}^{-n}(p_k+q_k)\leq \frac{r_{-1}}{2}.\end{equation}
Since  $p'_{-1}$, $p_{-i}$, $q_{-i}$  are integer we obtain that
  $p'_{-1}+\sum_{k=-2}^{-n}(p_k+q_k)$  takes the maximal value not  $\frac{r_{-1}}{2}$ but $\frac{r_{-1}}{2}-1$.

By Proposition \ref{p1} if one decomposes a vector corresponding to
a tableau with $m'_{-1,n}>0$ by the base
 \eqref{solindib0} then only monomials with  $p'_{-1}+\sum_{k=-2}^{-n}(p_k+q_k)\leq  m_{-1,n}-1=\frac{r_{-1}}{2}-1$ appear.
  Thus the image of \eqref{sop2} contains the span of vectors corresponding to tableaux with  $m'_{-1,n}>0$.

To prove the inverse embedding let us calculate the dimension of the
image and of the span.

Цe obtained that the  inequality \eqref{solindib0} for 
 $p_{-1}$,
$p_{k},q_{k}$ and $r_{-n},...,r_{-2},r_{-1}$ are equivalent to
inequalities \eqref{solg0} for $p'_{-1}$, $p_{k},q_{k}$ and
$r_{-n},...,r_{-2},\frac{r_{-1}}{2}-1$.

The corresponding monomials \eqref{solg0}, given by $p'_{-1}$, $p_{k},q_{k}$ and
$r_{-n},...,r_{-2},\frac{r_{-1}}{2}-1$,  define a base in the space
 of the problem of restriction   $\mathfrak{gl}_{n+1}\downarrow\mathfrak{gl}_{n-1}$ for the highest weight $[m_{-n}-1,...m_{-1}-1,0]$. In the same space  there exists another base
 indexed by tableau, which elements are integers

\begin{displaymath}
\label{gccgg}
 \xymatrix{ m_{-n}-1   && m_{-n+1}-1 &  ... & m_{-1}-1 && 0\\
 &\bar{m}'_{-n,n}   && ...   && \bar{m}'_{-1,n}\\
&& \bar{m}_{-n,n-1} &...& \bar{m}_{-2,n-1}\\
}
\end{displaymath}

Hence the dimension of the image of \eqref{sop2} is the number of
such tableaux.

To each such a tableaux there corresponds a tableau composed of
integers

\begin{displaymath}
 \xymatrix{ m_{-n}   && m_{-n+1} &  ... & m_{-1} && 0\\
 &\bar{m}'_{-n,n}+1   && ...   && \bar{m}'_{-1,n}+1\\
&& \bar{m}_{-n,n-1}+1 &...& \bar{m}_{-2,n-1}+1\\
}
\end{displaymath}

In this tableaux  $m'_{-1,n}=\bar{m}'_{-1,n}+1>0$ and each tableau и
with
 $m'_{-1,n}>0$  can be written in this manner. Hence the dimension of linear span of tableaux   with
   $m'_{-1,n}>0$ equals to the dimension of the image of the correspondence \eqref{sop2}.
   Thus they are equal
\endproof

\section{The Gelfand-Tsetlin-Zhelobenko base in the space of $\mathfrak{o}_{2n-1}$-highest vectors}

\label{gczs1}

In the previous Section we investigated the mapping

\begin{equation}
\label{sooti}
(p_{-1},p_{k},q_{k})\mapsto (p'_{-1}=[\frac{p_{-1}}{2}],p_{k},q_{k}, \sigma =0,1),\,\,\,k=-2,...,-n,\,\,\,\,,
\end{equation}
where $\sigma$ is a residue of the division of  $p_{-1}$ by  $2$,
with relates to the solution \eqref{solindib0} corresponding to $r_{-n},...,r_{-2},r_{-1}$ the solution
\eqref{solg0} corresponding to $r_{-n},...,r_{-2},[\frac{r_{-1}}{2}]$.

Compare formulas  \eqref{solindib0} and  \eqref{solg0}. We obtain the following statement.

\begin{prop}

The mapping inverse to \eqref{sooti} (non
everywhere-defined), can be written as follows

\begin{equation}
\label{soot}
f(z_{-n,-1},z_{-n,1},...,z_{-2,-1},z_{-2,1},z_{-1,1})\mapsto z_{0,1}^{\sigma}f(z_{-n,-1},z_{-n,1},...,z_{-2,-1},z_{-2,1},\frac{1}{2}z_{0,1}^2)
\end{equation}

\end{prop}

Let us construct Gelfand-Tsetlin tableaux and Gelfand-Tselin base for the restriction problem  $\mathfrak{o}_{2n+1}\downarrow\mathfrak{o}_{2n-1}$. For this pupose we construct a special base in the solution space of the restriction problem   $\mathfrak{gl}_{n+1}\downarrow\mathfrak{gl}_{n-1}$ that occur on the  side in  \eqref{soot}.

 Suppose that the highest weight as integer. Then in the space of solutions of the restriction problem   $\mathfrak{gl}_{n+1}\downarrow\mathfrak{gl}_{n-1}$ c $r_{-n},...,r_{-2},[\frac{r_{-1}}{2}]=m_{-1}$ there exist a base encoded by tableaux

\begin{displaymath}
 \xymatrix{ m_{-n}   && m_{-n+1} &  ... & m_{-1} && 0\\
 &m'_{-n,n}   && ...   && m'_{-1,n}\\
&& m_{-n,n-1} &...& m_{-2,n-1}\\
}
\end{displaymath}

Thus in this case in the solution space of the restriction problem  $\mathfrak{o}_{2n+1}\downarrow\mathfrak{o}_{2n-1}$ there exists a base encoded by such tableaux and a number $\sigma=0,1$.

Suppose that the highest weight is half integer. In the solution space of the restriction problem  $\mathfrak{gl}_{n+1}\downarrow\mathfrak{gl}_{n-1}$ c $r_{-n},...,r_{-2},[\frac{r_{-1}}{2}]=m_{-1}-1/2$ there exists a base encoded by such tableaux

\begin{displaymath}
 \xymatrix{ m_{-n}-1/2   && m_{-n+1}-1/2 &  ... & m_{-1}-1/2 && 0\\
 &\bar{m}'_{-n,n}   && ...   && \bar{m}'_{-1,n}\\
&& \bar{m}_{-n,n-1} &...& \bar{m}_{-2,n-1}\\
}
\end{displaymath}

Chose another indexation of this base. All element of this tableau is integer. Let us add to each element of this tableau   $1/2$ and obtain a tableau with half-integer elements.
In particular the lower row in a collection of eigenvalues of  $E_{-i,-i}+\frac{1}{2}id$, $i=2,...,n$.

 Thus in this case  in the solution space of the restriction problem  $\mathfrak{o}_{2n+1}\downarrow\mathfrak{o}_{2n-1}$  there exist a base encoded by tableau of the same type as in th case of integer highest weight but with half-integer elements and a $\sigma=0,1$. Thus we have proved the Theorem.

Proposition    \ref{pr5}  we described the domain of definition of  \eqref{soot}. Thus we come to the Theorem.

\begin{thm}

\label{osnt}

Let  $m_{-n,n}:=m_{-n},...,m_{-1,n}:=m_{-1}$.  Then in the space of
$\mathfrak{o}_{2n-1}$-highest vectors in a
$\mathfrak{o}_{2n+1}$-representation $V$  there exists a base indexed
by tableaux

\begin{align}\label{dia}\begin{split}&m_{-n,n}\geq m'_{-n,n}\geq m_{-n+1,n}\geq m'_{-n+1,n}\geq ...\geq m_{-1,n}\geq m'_{-1,n}\geq
0\\\sigma\,\,\, &m'_{-n,n}\geq m_{-n,n-1}\geq m'_{-n+1,n}\geq m_{-n+1,n-1}\geq
...\geq m_{-2,n-1}\geq m'_{-1,n}\end{split}.\end{align}

Here $\sigma$ takes only values $1$ and $0$, and other numbers are
simultaneously integers or half-integers. If the highest weight is integer and $m'_{-1,n}=0$ then
$\sigma=0$.
\end{thm}

As in the case $\mathfrak{sp}_{2n}$ the following statement takes
place.

\begin{prop}\label{strokaves1}
 The lower row of the tableau \eqref{dia} is a $\mathfrak{o}_{2n-1}$-weight of the corresponding $\mathfrak{o}_{2n-1}$-highest vector.
\end{prop}
\proof

In the case of integer highest weight the correspondence  \eqref{soot} conjugates the actions of  $E_{-i,-i}$  and $F_{-i,-i}$ for  $i=n,...,2$.  This follows form formulas  \eqref{fii} and \eqref{eii}.
  The lower row of a tableau appearing in the restriction problem  $\mathfrak{gl}_{n+1}\downarrow\mathfrak{gl}_{n-1}$ is a collection of eigenvalues of  $E_{-i,-i}$ for  $i=n,...,2$.
Hence after application of \eqref{soot} to this tableau we obtain a vector encoded by \eqref{dia} and the lower row of \eqref{dia}  is a collection of eigenvalues of  $F_{-i,-i}$ for  $i=n,...,2$.

In the case of half-integer highest weight the formulas \eqref{fii} and \eqref{eii} show that the correspondence \eqref{soot} conjugates the actions $F_{-i,-i}f$ and $E_{-i,-i}f+\frac{1}{2}f$. But in this case the tableau  \eqref{dia} is constructed in such way that  $m_{-i,n-1}$ is an eigenvalue of   $E_{-i,-i}f+\frac{1}{2}f$.

\endproof

Using Theorem \ref{osnt} and Proposition \ref{strokaves1} one can
construct the Gelfand-Tsetlin-Zhelobenko base vectors in a
representation of $\mathfrak{o}_{2n+1}$.


Let us write the formula for the $(-1)$-component of the weight of
the vector defined by a tableau
\begin{prop} The $(-1)$-component of the weight equals
\begin{equation}
\label{wess}
\sigma+2\sum_{k=1}^{n}m'_{-k,n}-\sum_{k=1}^{n}m_{-k,n}+\sum_{k=2}^n m_{-k,n-1}.
\end{equation}
\end{prop}

\proof
Let the highest weight be integer. By \eqref{fii} and \eqref{eii} one has

\begin{align*}
&F_{-1,-1}f=\sum_{i=2}^n(z_{-i,-1}\frac{\partial}{\partial z_{-i,-1}}-z_{-i,1}\frac{\partial}{\partial z_{-i,1}})f-z_{0,1}\frac{\partial}{\partial z_{0,1}}f-m_{-1}f,\\
&(E_{-1,-1}-E_{1,1})f=\sum_{i=2}^n(z_{-i,-1}\frac{\partial}{\partial z_{-i,-1}}-z_{-i,1}\frac{\partial}{\partial z_{-i,1}})f-2z_{-1,1}\frac{\partial}{\partial z_{-1,1}}f-m_{-1}f.
\end{align*}

Let $\sigma=0$. The under the correspondence \eqref{soot} the change of variables $z_{-1,1}\mapsto \frac{z_{0,1}^2}{2}$ is performed.  Thus the actions  $2z_{-1,1}\frac{\partial}{\partial z_{-1,1}}$ and  $z_{0,1}\frac{\partial}{\partial z_{0,1}}$  are conjugated. Hence the actions of  $F_{-1,-1}$ and $E_{-1,-1}-E_{1,1}$ are conjugated.
 Thus the eigenvalue of $F_{-1,-1}$ on the vector  \eqref{dia} is a difference of eigenvalues of $E_{-1,-1}$ and  $E_{1,1}$.
These eigenvalues are equal to $\sum_{k=1}^{n}m'_{-k,n}-\sum_{k=2}^n m_{-k,n-1}$ and $\sum_{k=1}^{n}m_{-k,n}-\sum_{k=1}^{n}m'_{-k,n}$. The difference of these expressions is \eqref{wess}.

In the case  $\sigma=1$ при соответствии \eqref{soot} происходит еще умножение на  $z_{0,1}$. Поэтому к разности собственных значений $E_{-1,-1}$ и  $E_{1,1}$ добавляется  $1$.

Suppose that the weight is half-integer. Let  $\sigma=0$.  Then by formulas \eqref{fii} and \eqref{eii} under the action of  \eqref{soot}  the actions of  $F_{-1,-1}$ and $E_{-1,-1}-E_{1,1}+\frac{1}{2}id$ are conjugated. In term of the tableau  \eqref{dia} the eigenvalues of $E_{-1,-1}$ and $E_{1,1}$ are equal   $\sum_{k=1}^{n}m'_{-k,n}-\sum_{k=2}^n m_{-k,n-1}-\frac{1}{2}$ and $\sum_{k=1}^{n}m_{-k,n}-\sum_{k=1}^{n}m'_{-k,n}$. Their difference plus $\frac{1}{2}$ equals to \eqref{wess}.

In the case  $\sigma=1$ to this expression $1$ must be added.

\endproof


\end{document}